\documentclass[12pt,reqno,openbib,runningheads,a4paper]{amsart}
\usepackage{graphicx}
\usepackage{amssymb}
\usepackage{amsfonts}
\usepackage{amsmath}
\usepackage{graphicx}
\usepackage[authoryear]{natbib}
\usepackage{epsf}
\usepackage{epstopdf}
\usepackage[legalpaper,bookmarks=true,colorlinks=true,linkcolor=blue,citecolor=blue]{hyperref}

\providecommand{\U}[1]{\protect\rule{.1in}{.1in}}
\newtheorem{proposition}{Proposition}[section]

\newtheorem{remark}{Remark}[section]

\makeatletter
\renewcommand{\@biblabel}[1]{}
\makeatother
\input{tcilatex}

\begin{document}

\begin{center}
{\Large \textbf{Involving copula functions in Conditional Tail Expectation}}%
\medskip \medskip

{\large Brahim Brahimi}\footnote{%
\noindent E-mail addresses:{\small {\texttt{brah.brahim@gmail.com},
Tel.:+213-7 73 54 60 63; fax:+213-33 74 77 88.}}}\medskip

{\small \textit{Laboratory of Applied Mathematics, Mohamed Khider
University, Biskra, Algeria}}\bigskip \medskip
\end{center}

\noindent\textbf{Abstract}\medskip

\noindent Our goal in this paper is to propose an alternative risk measure
which takes into account the fluctuations of losses and possible
correlations between random variables. This new notion of risk measures,
that we call Copula Conditional Tail Expectation describes the expected
amount of risk that can be experienced given that a potential bivariate risk
exceeds a bivariate threshold value, and provides an important measure for
right-tail risk. An application to real financial data is given.\medskip


\section{\textbf{Introduction}}

\noindent In actuarial science literature a several risk measures have been
proposed, namely: the Value-at-Risk (VaR), the expected shortfall or the
conditional tail expectation (CTE), the distorted risk measures (DRM), and
recently the copula distorted risk measure (CDRM) as risk measure which
takes into account the fluctuations dependence between random variables
(rv). See \cite{BrMN10}.$\medskip $

\noindent The CTE in risk analysis represents the conditional expected loss
given that the loss exceeds its VaR and provides an important measure for
right-tail risk. In this paper we will always consider random variables with
finite mean. For a real number $s\ $in $\left( 0,1\right) ,$ the CTE of a
risk $X$ is given by%
\begin{equation}
\mathbb{CTE}\left( s\right) :=\mathbb{E}\left[ \left. X\right\vert
X>VaR_{X}\left( s\right) \right] ,  \label{e1}
\end{equation}
where $VaR_{X}\left( s\right) :=\inf\left\{ x:F\left( x\right) \geq
s\right\} $ is the quantile of order $s$ pertaining to distribution function
(df) $F.$\medskip

\noindent One of the strategy of an Insurance companies is to set aside
amounts of capital from which it can draw from in the event that premium
revenues become insufficient to pay out claims. Of course, determining these
amounts is not a simple calculation. It has to determine the best risk
measure that can be used to determine the amount of loss to cover with a
high degree of confidence.$\medskip$

\noindent Suppose now that we deal with a couple of random losses $%
(X_{1},X_{2}).$ It's clear that the CTE of $X_{1}$ is unrelated with $X_{2}.$
If we had to control the overflow of the two risks $X_{1}$ and $X_{2}$ at
the same time, CTE does not answer the problem, then we need another
formulation of CTE which takes into account the excess of the two risks $%
X_{1}$ and $X_{2}.$ Then we deal with the amount
\begin{equation}
\mathbb{E}\left[ \left. X_{1}\right\vert X_{1}>VaR_{X_{1}}\left( s\right)
,X_{2}>VaR_{X_{2}}\left( t\right) \right] .  \label{e2}
\end{equation}%
If the couple of random losses $(X_{1},X_{2})$ are independents rv's then
the amount (\ref{e2}) defined only the CTE of $X_{1}.$ Therefore the case of
independence is not important.\medskip

\noindent In the recent years dependence is beginning to play an important
role in the world of risk. The increasing complexity of insurance and
financial activities products has led to increased actuarial and financial
interest in the modeling of dependent risks. While independence can be
defined in only one way, but dependence can be formulated in an infinite
ways. Therefore, the assumption of independence it makes the treatment easy.
Nevertheless, in applications dependence is the rule and independence is the
exception.\medskip

\noindent The copulas is a function that completely describes the dependence
structure, it contains all the information to link the marginal
distributions to their joint distribution. To obtain a valid multivariate
distribution function, we combines several marginal distribution functions,
or a different distributional families, with any copula function. Using
Sklar's theorem \citep{Sklar59}, we can construct a bivariate distributions
with arbitrary marginal distributions. Thus, for the purposes of statistical
modeling, it is desirable to have a large collection of copulas at one's
disposal. A great many examples of copulas can be found in the literature,
most are members of families with one or more real parameters. For a formal
treatment of copulas and their properties, see the monographs by \cite%
{HutLai90}, \cite{Dall91}, \cite{Joe97}, the conference proceedings edited
by \cite{BeS97}, \cite{Cua02}, \cite{Dha03} and the textbook of \cite{Nels06}%
.\medskip

\noindent Recently in finance, insurance and risk management has emphasized
the importance of positive or negative quadrant dependence notions (PQD or
NQD) introduced by \cite{LEH66}, in different areas of applied probability
and statistics, as an example, see; \cite{Dhaen97}, \cite{Denuitetal01}. Two
rv's are said to be PQD when the probability that they are simultaneously
large (or small) is at least as great as it would be were they are
independent. In terms of copula, if their copula is greater than their
product, i.e., $C(u_{1},u_{2})>u_{1}u_{2}$ or, simply $C>C^{\perp },$ where $%
C^{\perp }$ denotes the product copula. For the sake of brevity, we will
restrict ourselves to concepts of positive dependence. \medskip

\noindent The main idea of this paper is to use the information of
dependence between PQD or NQD risks to quantifying insurance losses and
measuring financial risk assessments, we propose a risk measure defined by:%
\begin{equation*}
\mathbb{CCTE}_{X_{1}}\left( s;t\right) :=\mathbb{E}\left[ \left.
X_{1}\right\vert X_{1}>VaR_{X_{1}}\left( s\right) ,X_{2}>VaR_{X_{2}}\left(
t\right) \right] .
\end{equation*}%
We will call this new risk measure by the \textit{Copula Conditional Tail
Expectation} (CCTE), like a risk measure which measure the conditional
expectation given the two dependents losses exceeds $VaR_{X_{1}}\left(
s\right) $ and $VaR_{X_{2}}\left( t\right) $ for $0<s,t<1$ and usually with $%
s,t>0.9.$ Again, CCTE satisfies all the desirable properties of a coherent
risk measure \citep{ARZ99}. The notion of copula in risk measure filed has
recently been considered by several authors, see for instance \cite{Emb03},
\cite{ANN04}, \cite{Luci09}, \cite{BrMN10} and the references
therein.\medskip

\noindent This risk measures can give a good quantifying of losses when we
have a combined dependents risks, this dependence can influence in the
losses of interested risks. Therefore, quantify the riskiness of our
position is useful to decide if it acceptable or not. For this reason we use
the all informations a bout this interest risk. The dependence of our risk
with other risks is one of important information that we must take it in
consideration.\medskip

\noindent This paper is organized as follows. In section \ref{section2}, we
give an explicit formulations of the new notion copula conditional tail
expectation risk measure in bivariate case. The relationship of this new
concept and tail dependence measure, given in section \ref{section3}. In
section \ref{Section example} we presents an illustration examples to
explain how to use the new CCTE measure. Application in real financial data
is given in section \ref{Applicaltion}. Concluding notes are given in
Section \ref{Conclu}. Proofs are relegated to the Appendix.

\section{\label{section2}\textbf{Copula conditional tail expectation}}

\noindent A risk measure quantifies the risk exposure in a way that is
meaningful for the problem at hand. The most commonly used risk measure in
finance and insurance are: VaR and CTE (also known as Tail-VaR or expected
shortfall). The risk measure is simply the loss size for which there is a
small (e.g. $1\%$) probability of exceeding. For some time, it has been
recognized that this measure suffers from serious deficiencies if losses are
not normally distributed.\medskip

\noindent According to \cite{ARZ99} and \cite{WirchHard99}, the conditional
tail expectation of a random variable $X_{1}$ at its $VaR_{X_{1}}\left(
s\right) $ is defined by:%
\begin{equation*}
\mathbb{CTE}_{X_{1}}\left( s\right) =\frac{1}{1-F_{X_{1}}(VaR_{X_{1}}\left(
s\right) )}\int_{VaR_{X_{1}}\left( s\right) }^{\infty }xdF_{X_{1}}(x),
\end{equation*}%
where $F_{X_{1}}$ is the df of $X_{1}$.\medskip

\noindent Since $X_{1}$ is continuous, then $F_{X_{1}}(VaR_{X_{1}}\left(
s\right) )=s,$ it follows that for all $0<s<1$
\begin{equation}
\mathbb{CTE}_{X_{1}}\left( s\right) =\frac{1}{1-s}\int_{s}^{1}VaR_{X_{1}}%
\left( u\right) du.  \label{eq2}
\end{equation}%
The CTE can be larger that the VaR measure for the same value of level $s$
described above since it can be thought of as the sum of the quantile $%
VaR_{X_{1}}\left( s\right) $ and the expected excess loss. Tail-VaR is a
coherent measure in the sense of \cite{ARZ99}. For the application of this
kind of coherent risk measures we refer to the papers \cite{ARZ99} and \cite%
{WirchHard99}.\medskip

\noindent Thus the CTE is nothing, see \cite{OverSko08}, but the
mathematical transcription of the concept of "average loss in the worst $%
100(1-s)\%$ case", defining by $\upsilon =VaR_{X_{1}}(s)$ a critical loss
threshold corresponding to some confidence level $s,$ $\mathbb{CTE}%
_{X_{1}}(s)$ provides a cushion against the mean value of losses exceeding
the critical threshold $\upsilon .$\medskip

\noindent Now, assume that $X_{1}$ and $X_{2}$ are dependent with joint df $%
H $ and continuous margins $F_{X_{i}},$ $i=1,2,$ respectively. Through this
paper we calls $X_{1}$ the \textit{target risk} and $X_{2}$ the \textit{%
associated risk}. In this case, the problem becomes different and its
resolution requires more than the usual background.\medskip

\noindent Our contribution is to introduce the copula notion to provide more
flexibility to the CTE of risk of rv's in terms of loss and dependence
structure. For comprehensive details on copulas one may consult the textbook
of \cite{Nels06}.\medskip

\noindent According to Sklar's Theorem \cite{Sklar59}, there exists a unique
copula $C:\left[ 0,1\right] ^{d}\rightarrow \left[ 0,1\right] $ such that
\begin{equation}
H\left( x_{1},x_{2}\right) =C\left( F_{1}\left( x_{1}\right) ,F_{2}\left(
x_{2}\right) \right) .  \label{sklar}
\end{equation}%
The formula of CTE focuses only on the average of loss. For this you should
think of an other formula to be more inclusive, this formula must take in
consideration the dependence structure and the behavior of margin tails.
These two aspects have an important influence when quantifying risks. On the
other hand if the correlation factor is neglected, the calculation of the
CTE follows from formula (\ref{eq2}), which only focuses on the target
risk.\medskip

\noindent Now by considering the correlation between the target and the
associated risks, we define a new notion of CTE called \textit{Copula
Conditional Tail Expectation} (CCTE) given in (\ref{e2}), this notion led to
give a risk measurement focused in the target risk and the link between
target and associated risk.\medskip

\noindent Let's denote the survival functions by $\overline{F}_{i}(x_{i})=%
\mathbb{P}(X_{i}>x_{i}),$ $i=1,2,$ and the joint survival function by $%
\overline{H}(x_{1},x_{2})=\mathbb{P}(X_{1}>x_{1},X_{2}>x_{2}).$ The function
$\overline{C}$ which couples $\overline{H}$ to $\overline{F}_{i},$ $i=1,2$
via $\overline{H}(x_{1},x_{2})=\overline{C}(\overline{F}_{1}(x_{1}),%
\overline{F}_{2}(x_{2}))$ is called the survival copula of $\left(
X_{1},X_{2}\right) .$ Furthermore, $\overline{C}$ is a copula, and
\begin{equation}
\overline{C}(u_{1},u_{2})=u_{1}+u_{2}-1+C(1-u_{1},1-u_{2}),  \label{SC}
\end{equation}%
where $C$ is the (ordinary) copula of $X_{1}$ and $X_{2}.$ For more details
on the survival copula function see, Section 2.6 in \cite{Nels06}.\medskip

\noindent The CCTE of the target risk $X_{1}$ with respect to the associated
risk $X_{2}$ is given in the following proposition.\bigskip

\begin{proposition}
\label{Prop1}Let $\left( X_{1},X_{2}\right) $ a bivariate rv with joint df
represented by the copula $C.$ Assume that $X_{1}$ have a finite mean and df
$F_{X_{1}}.$ Then for all $s$ and $t$ in $\left( 0,1\right) $ the copula
conditional tail expected of $X_{1}$ with respect to the bivariate
thresholds $(s,t)$ is given by%
\begin{equation}
\mathbb{CCTE}_{X_{1}}\left( s;t\right) =\frac{\int_{s}^{1}F_{X_{1}}^{-1}%
\left( u\right) \left( 1-C_{u}(u,t)\right) du}{\overline{C}\left(
1-s,1-t\right) },  \label{CCTE}
\end{equation}%
where $F_{X_{1}}^{-1}$ is the quantile function of $F_{X_{1}}$and $%
C_{u}(u,v):=\partial C(u,v)/\partial u.$
\end{proposition}

\noindent This notion does not depend on the df of the associated risk, but
it depend only by the copula function and the df of target risk.\medskip

\noindent By definition of PQD risks we have that $C(u,v)>uv,$ then it easy
to check that
\begin{equation*}
\mathbb{CCTE}_{X_{1}}\left( s;t\right) \leq \frac{\mathbb{CTE}_{X_{1}}\left(
s\right) }{\left( 1-t\right) }\text{ for }s,t<1,
\end{equation*}%
next, in Section \ref{Section example}, we will proved that the risk when we
consider the correlation between PQD risks is greater than in the case of a
single one. That means, for all $s\leq t$ and $s,$ $t$ in $\left( 0,1\right)
$ then
\begin{equation}
\mathbb{CCTE}_{X_{1}}\left( s;t\right) \geq \mathbb{CTE}_{X_{1}}\left(
s\right) .  \label{CCTE theorem}
\end{equation}%
Notice that in the NQD rv's we have the reverse inequality of (\ref{CCTE
theorem}) and the CCTE coincide with CTE measures in the non-dependence
case, i.e. the copula $C=C^{\bot }.$

\section{CCTE and tail dependence\label{section3}}

\noindent The concept of tail dependence is an asymptotic measure of the
dependence between two random variables in the tail of their joint
distribution function. Specifically, tail dependence is the probability that
a random variable $X_{1}$ and $X_{2}$ takes a values in the extreme tail of
its distributions simultaneously, for example we consider $X_{1}$ and $X_{2}$
which measure bankruptcy for two companies and both companies simultaneously
go bankrupt.\medskip

\noindent We describes the joint upper tail dependence of the random
variables $X_{1}$ and $X_{2}:$%
\begin{equation*}
\lim_{\substack{ t\rightarrow 1^{-}  \\ s\rightarrow 1^{-}}}\mathbb{P}\left(
\left. X_{1}>F_{X_{1}}^{-1}\left( s\right) \right\vert
X_{2}>F_{X_{2}}^{-1}\left( t\right) \right)
\end{equation*}%
However, it can be seen as a good indicator of systemic risk (for $s=t$). If
we considering the tail dependence as a dependence measure in the extreme
tails of the joint distribution, it is possible for two rv's to be
dependent, but for there to be no dependence in the tail of the
distributions, this is the case described for example by the Gaussian
copula, hyperbolic copula or Farlie-Gumbel-Morgenstern copula (tail
dependence is zero). Furthermore, the Clayton copula puts the entire tail
dependence in the lower tail unlike Gumbel copula in the upper tail and the
Student copula behave identically in the lower as in the upper tail.
However, it is not suitable to model extreme negative outcomes similarly as
with extreme positive outcomes.\medskip

\begin{remark}
Negative outcomes can be treated in the same way that the extremes positive
outcomes by replacing their copula by the survival copula.
\end{remark}

\noindent The tail dependence can be also expressed through copula%
\begin{equation*}
\lambda _{U}=\lim_{u\rightarrow 1^{-}}\frac{1-2u+C\left( u,u\right) }{1-u}%
\text{ and }\lambda _{L}=\lim_{u\rightarrow 0^{+}}\frac{C\left( u,u\right) }{%
u}.
\end{equation*}%
Now, let's denote by
\begin{equation*}
\tilde{\lambda}_{U}\left( u,v\right) :=\frac{1-u-v+C\left( u,v\right) }{1-v}%
\text{ and }\tilde{\lambda}_{L}\left( u,v\right) :=\frac{C\left( u,v\right)
}{v}.
\end{equation*}%
Note that $\lim_{u,v\rightarrow 1^{-}}\tilde{\lambda}_{U}\left( u,v\right)
=\lambda _{U}$ and $\lim_{u,v\rightarrow 0^{+}}\tilde{\lambda}_{L}\left(
u,v\right) =\lambda _{L}.$ We can rewrite CCTE of according to $\tilde{%
\lambda}_{U}$ as%
\begin{equation*}
\mathbb{CCTE}_{X_{1}}\left( s;t\right) =\frac{\int_{s}^{1}F_{X_{1}}^{-1}%
\left( u\right) \left( 1-C_{u}(u,t)\right) du}{\left( 1-t\right) \tilde{%
\lambda}_{U}\left( s,t\right) },
\end{equation*}%
this has no impact on the limiting value at $0$ for PQD risks. Then we have
\begin{equation*}
\lim_{\substack{ s\rightarrow 1^{-}  \\ t\rightarrow 1^{-}}}\left(
1-t\right) \tilde{\lambda}_{U}\left( s,t\right) =0.
\end{equation*}%
From Theorem 2.2.7 in \cite[page 13]{Nels06} we have $0\leq C_{u}(u,t)\leq 1$
for such $u$ and $t,$ then%
\begin{equation*}
\left\vert \mathbb{CCTE}_{X_{1}}\left( s;t\right) \right\vert \leq
\left\vert \frac{\int_{s}^{1}F_{X_{1}}^{-1}\left( u\right) du}{\left(
1-t\right) \tilde{\lambda}_{U}\left( s,t\right) }\right\vert \leq \left\vert
\frac{\mathbb{E}\left( X_{1}\right) }{\left( 1-t\right) \tilde{\lambda}%
_{U}\left( s,t\right) }\right\vert .
\end{equation*}%

\noindent In the next section, we give an example to describe the impact of
the upper tail dependence nearly $1$ and the lower tail dependence near $0$
in CCTE, and we discuss the relationship between the properties of the
dependence of copula model with upper and lower tail dependence and how to
derive the degree of risk in each case. \medskip

\section{\textbf{Illustration examples}\label{Section example}}

\subsection{CCTE via Farlie-Gumbel-Morgenstern Copulas\label{FGMsubsection}}

\noindent One of the most important parametric family of copulas is the
Farlie-Gumbel-Morgenstern (FGM) family defined as
\begin{equation}
C_{\theta }^{FGM}(u,v)=uv+\theta uv(1-u)(1-v),\ \ \ u,v\in \lbrack 0,1],
\label{FGM}
\end{equation}%
where $\theta \in \lbrack -1,1].$ The family was discussed by \cite{Mor56},
\cite{Gumb58} and \cite{Farlie60}.\medskip

\noindent The copula given in (\ref{FGM}) is PQD for $\theta \in (0,1]$ and
NQD for $\theta \in \lbrack -1,0).$ In practical applications this copula
has been shown to be somewhat limited, for copula dependence parameter $%
\theta \in \left[ -1,1\right] ,$ Spearman's correlation $\rho \in \left[
-1/3,1/3\right] $ and Kendall's $\tau \in \left[ -2/9,2/9\right] ,$ for more
details on copulas see, for example, \cite{Nels06}.\medskip

\noindent Members of the FGM family are symmetric, i.e., $%
C_{\theta}^{FGM}(u,v)=C_{\theta}^{FGM}(v,u)$ for all $(u,v)$ in $\left[ 0,1%
\right] ^{2}$ and have the lower and upper tail dependence coefficients
equal to $0.$\medskip

\noindent A pair $\left( X,Y\right) $ of rv's is said to be exchangeable if
the vectors $\left( X,Y\right) $ and $\left( Y,X\right) $ are identically
distributed. Note that, in applications, exchangeability may not always be a
realistic assumption. For identically distributed continuous random
variables, exchangeability is equivalent to the symmetry of the FGM
copula.\medskip\

\noindent For practical purposes we consider a copula families with only
positive dependence. Furthermore, risk models are often designed to model
positive dependence, since in some sense it is the \textquotedblleft
dangerous\textquotedblright\ dependence: assets (or risks) move in the same
direction in periods of extreme events, see \cite{EmLMC03}. \medskip

\noindent Consider the bivariate loss PQD rv's $(X_{i},Y),$ $i=1,2,3,$
having continuous marginal df's $F_{X_{i}}(x)$ and $F_{Y}(y)$ and joint df $%
H_{X_{i},Y}(x,y) $ represented by FGM copula of parameters $\theta_{i}$,
respectively for $i=1,2,3$%
\begin{equation*}
H_{X_{i},Y}(x,y)=C_{\theta_{i}}^{FGM}(F_{X_{i}}\left( x\right) ,F_{Y}\left(
y\right) ).
\end{equation*}
The marginal survival functions $\overline{F}_{X_{i}}(x),$ $i=1,2,3$ and $%
\overline{F}_{Y}(y)$ are given by
\begin{equation}
\overline{F}_{X_{i}}\left( x\right) =\left\{
\begin{tabular}{ll}
$\left( 1+x\right) ^{-\alpha},$ & $x\geq0,$ \\
$1,$ & $x<0,$%
\end{tabular}
\ \ \ \right. \text{ and }\ \ \overline{F}_{Y}\left( y\right) =\left\{
\begin{tabular}{ll}
$\left( 1+y\right) ^{-\alpha},$ & $y\geq0,$ \\
$1,$ & $y<0.$%
\end{tabular}
\ \ \ \ \ \ \ \right.  \label{F}
\end{equation}
where $\alpha>0$ called the Pareto index, the case $\alpha\in(1,2)$ means
that $X_{i}$ have a heavy-tailed distributions. So that $X_{i}$ and $Y$ have
identical Pareto df's.\medskip

\noindent For each couple $\left( X_{i},Y\right) ,$ $i=1,2,3,$ we propose $%
\theta_{1}=0.01,$ $\theta_{2}=0.5$ and $\theta_{3}=1,$ respectively. The
choice of parameters $\theta_{i},i=1,2,3$ correspond respectively to the
weak, medium and the high dependence.\medskip

\noindent In this example, the target risks are $X_{i}$ and the associated
risk is $Y.$ The $\mathbb{CTE}$'s and the VaR's of $X_{i}$ are the same and
are given respectively by%
\begin{equation}
\mathbb{CTE}_{X_{i}}\left( s\right) =\frac{\alpha\left( 1-s\right)
^{-1/\alpha}}{\alpha-1}  \label{cteXi}
\end{equation}
and
\begin{equation}
VaR_{X_{i}}\left( s\right) =(1-s)^{-1/\alpha},  \label{VaR}
\end{equation}
for $i=1,2,3.$\medskip

\noindent We have that%
\begin{equation}
\overline{C}\left( 1-s,1-t\right) =(1-s)(1-t)\left( st\theta _{i}+1\right) .
\label{aa}
\end{equation}%
Now, we calculate%
\begin{equation*}
\mathbb{CCTE}_{X_{1}}\left( s;t\right) =\frac{1}{\overline{C}\left(
1-s,1-t\right) }\int_{s}^{1}(1-u)^{-1/\alpha }\left( 1-t\right) \left(
2tu\theta _{i}-t\theta _{i}+1\right) du
\end{equation*}%
by substitution (\ref{aa}) we get
\begin{equation}
\mathbb{CCTE}_{X_{i}}\left( s;t\right) =\frac{\alpha \left( 2\alpha +t\theta
_{i}-2st\theta _{i}+2st\alpha \theta _{i}-1\right) }{\left( 2\alpha
^{2}-3\alpha +1\right) \left( st\theta _{i}+1\right) }\left( 1-s\right)
^{-1/\alpha }.  \label{FGMPro}
\end{equation}

\begin{figure}[h]
\begin{center}
\includegraphics[height=3.5243in,width=4.754in]{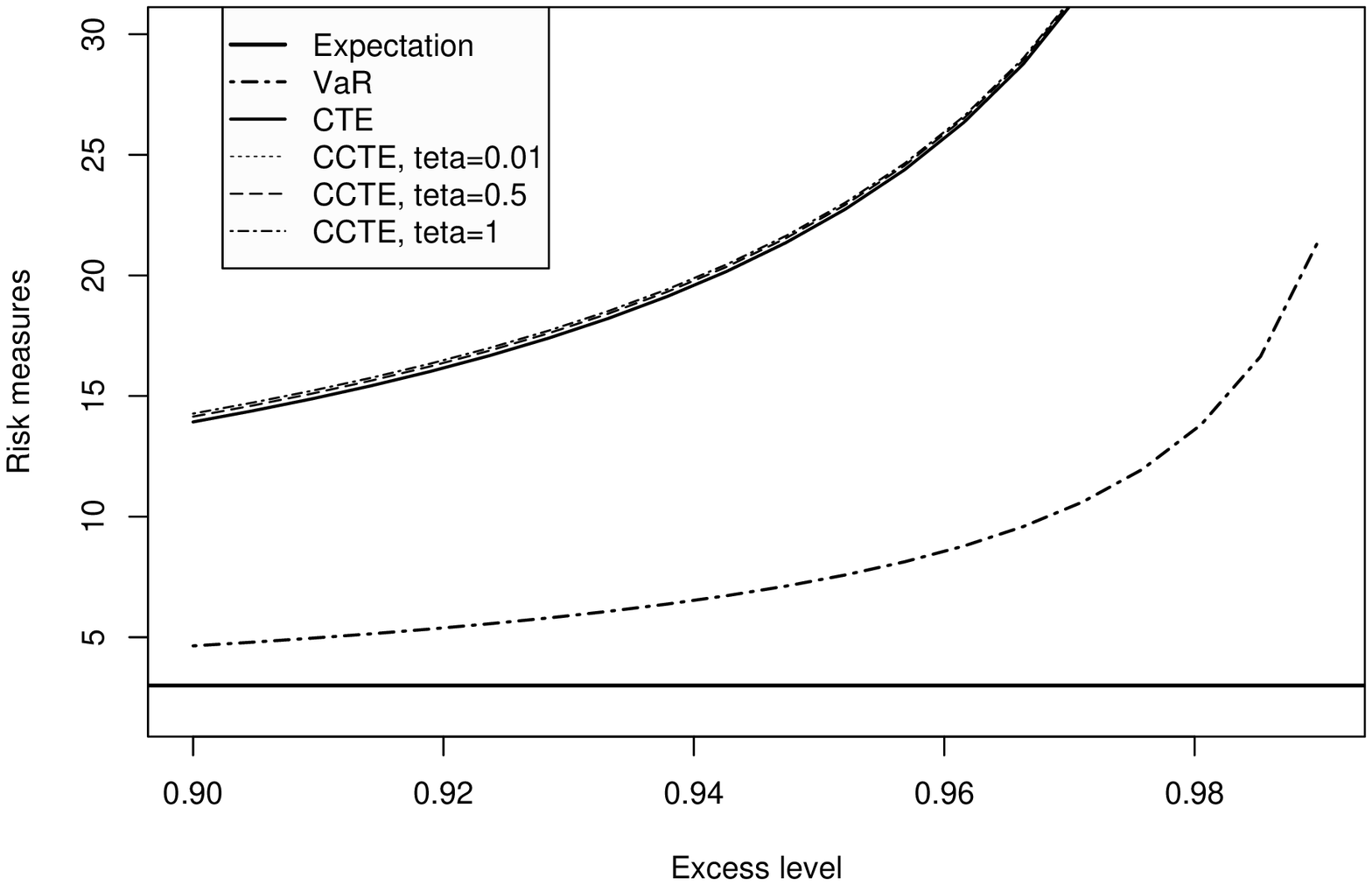}
\end{center}
\caption{$\mathbb{CCTE}$, $\mathbb{CTE}$ and $VaR$ risks measures of PQD
pareto (1.5) rv's with FGM copula and $0.9\leq s=t\leq0.99$}
\label{FigFGM}
\end{figure}

\noindent We have in Table \ref{TAB1} and Figures \ref{FigFGM}\ the
comparison of the riskiness of $X_{1},$ $X_{2}$ and $X_{3}.$ Recall that,
the $\mathbb{CTE}$'s risk measure of $X_{i}$ at level $s$ are the same in
all cases. Note that $\mathbb{CCTE}$ coincide with $\mathbb{CTE}$ in the
independence case $(\theta _{1}=0).$ The $\mathbb{CCTE}$ of the loss $X_{3}$
is riskier than $X_{2}$ and $X_{1}$ but not very significant, in the 6th
column of Table \ref{TAB1}, the relative difference between $64.7946$ and $%
64.633$ is only about $0.025\%$. This is due to that FGM copula does not
take into account the dependence in upper and lower tail $(\lambda
_{L}=\lambda _{U}=0).$ In this case we can not clearly confirm which is the
risk the more dangerous.

%TCIMACRO{\TeXButton{B}{\begin{table}[h] \centering} }%
%BeginExpansion
\begin{table}[h] \centering
%EndExpansion
\begin{tabular}{rccccc}
${\small s}$ & \multicolumn{1}{|c}{${\small 0.9000}$} & ${\small 0.9225}$ & $%
{\small 0.9450}$ & ${\small 0.9675}$ & ${\small 0.9900}$ \\ \hline
${\small VaR}_{X_{i}}\left( s\right) $ & \multicolumn{1}{|c}{${\small 4.6415}
$} & ${\small 5.5013}$ & ${\small 6.9144}$ & ${\small 9.8192}$ & ${\small %
21.5443}$ \\
$\mathbb{CTE}_{X_{i}}\left( s\right) $ & \multicolumn{1}{|c}{${\small 13.9247%
}$} & ${\small 16.5039}$ & ${\small 20.7433}$ & ${\small 29.4577}$ & $%
{\small 64.6330}$ \\ \hline
\multicolumn{1}{l}{${\small t}$} & \multicolumn{5}{|c}{$\mathbb{CCTE}%
_{X_{1}}\left( s,t\right) ,\ \ {\small \theta =0.01}$} \\ \hline
\multicolumn{1}{l}{${\small 0.9000}$} & \multicolumn{1}{|c}{${\small 13.9309}
$} & ${\small 13.9311}$ & ${\small 13.9312}$ & ${\small 13.9314}$ & ${\small %
13.9316}$ \\
\multicolumn{1}{l}{${\small 0.9225}$} & \multicolumn{1}{|c}{${\small 16.5096}
$} & ${\small 16.5097}$ & ${\small 16.5099}$ & ${\small 16.5100}$ & ${\small %
16.5101}$ \\
\multicolumn{1}{l}{${\small 0.9450}$} & \multicolumn{1}{|c}{${\small 20.7484}
$} & ${\small 20.7485}$ & ${\small 20.7487}$ & ${\small 20.7488}$ & ${\small %
20.7489}$ \\
\multicolumn{1}{l}{${\small 0.9675}$} & \multicolumn{1}{|c}{${\small 29.4619}
$} & ${\small 29.4620}$ & ${\small 29.4621}$ & ${\small 29.4623}$ & ${\small %
29.4624}$ \\
\multicolumn{1}{l}{${\small 0.9900}$} & \multicolumn{1}{|c}{${\small 64.6359}
$} & ${\small 64.6359}$ & ${\small 64.6360}$ & ${\small 64.6361}$ & ${\small %
64.6362}$ \\ \hline\hline
\multicolumn{1}{l}{${\small t}$} & \multicolumn{5}{|c}{$\mathbb{CCTE}%
_{X_{2}}\left( s,t\right) ,\ \ {\small \theta =0.5}$} \\ \hline
\multicolumn{1}{l}{${\small 0.9000}$} & \multicolumn{1}{|c}{${\small 14.1477}
$} & ${\small 14.1517}$ & ${\small 14.1555}$ & ${\small 14.1594}$ & ${\small %
14.1631}$ \\
\multicolumn{1}{l}{${\small 0.9225}$} & \multicolumn{1}{|c}{${\small 16.7072}
$} & ${\small 16.7108}$ & ${\small 16.7143}$ & ${\small 16.7178}$ & ${\small %
16.7212}$ \\
\multicolumn{1}{l}{${\small 0.9450}$} & \multicolumn{1}{|c}{${\small 20.9234}
$} & ${\small 20.9266}$ & ${\small 20.9297}$ & ${\small 20.9327}$ & ${\small %
20.9357}$ \\
\multicolumn{1}{l}{${\small 0.9675}$} & \multicolumn{1}{|c}{${\small 29.6077}
$} & ${\small 29.6103}$ & ${\small 29.6129}$ & ${\small 29.6154}$ & ${\small %
29.6179}$ \\
\multicolumn{1}{l}{${\small 0.9900}$} & \multicolumn{1}{|c}{${\small 64.7336}
$} & ${\small 64.7353}$ & ${\small 64.7370}$ & ${\small 64.7387}$ & ${\small %
64.7404}$ \\ \hline\hline
\multicolumn{1}{l}{${\small t}$} & \multicolumn{5}{|c}{$\mathbb{CCTE}%
_{X_{3}}\left( s,t\right) ,\ \ {\small \theta =1}$} \\ \hline
\multicolumn{1}{l}{${\small 0.9000}$} & \multicolumn{1}{|c}{${\small 14.2709}
$} & ${\small 14.2756}$ & ${\small 14.2803}$ & ${\small 14.2848}$ & ${\small %
14.2892}$ \\
\multicolumn{1}{l}{${\small 0.9225}$} & \multicolumn{1}{|c}{${\small 16.8183}
$} & ${\small 16.8226}$ & ${\small 16.8267}$ & ${\small 16.8308}$ & ${\small %
16.8348}$ \\
\multicolumn{1}{l}{${\small 0.9450}$} & \multicolumn{1}{|c}{${\small 21.0208}
$} & ${\small 21.0245}$ & ${\small 21.0281}$ & ${\small 21.0316}$ & ${\small %
21.0351}$ \\
\multicolumn{1}{l}{${\small 0.9675}$} & \multicolumn{1}{|c}{${\small 29.6880}
$} & ${\small 29.6910}$ & ${\small 29.6940}$ & ${\small 29.6969}$ & ${\small %
29.6997}$ \\
\multicolumn{1}{l}{${\small 0.9900}$} & \multicolumn{1}{|c}{${\small 64.7868}
$} & ${\small 64.7888}$ & ${\small 64.7908}$ & ${\small 64.7927}$ & ${\small %
64.7946}$ \\ \hline\hline
\multicolumn{1}{l}{} &  &  &  &  &
\end{tabular}%
\caption{Risk measures of dependent pareto(1.5) rv's with FGM copula.}\label{TAB1}%
\end{table}%

\subsection{CCTE via Archimedean Copulas}

\noindent A bivariate copula is said to be Archimedean \citep[see,][]{Gen86}
if it can be expressed by
\begin{equation*}
C(u_{1},u_{2})=\psi ^{\lbrack -1]}\left( \psi (u_{1})+\psi (u_{2})\right) ,
\end{equation*}%
where $\psi ,$ called the generator of $C,$ is a continuous strictly
decreasing convex function from $\left[ 0,1\right] $ to $[0,\infty ]$ such
that $\psi (1)=0$ with $\psi ^{\lbrack -1]}$ denotes the \textit{%
pseudo-inverse}\ of $\psi ,$ that is%
\begin{equation*}
\psi ^{\lbrack -1]}\left( t\right) =\left\{
\begin{tabular}{lll}
$\psi ^{-1}\left( t\right) ,$ & for & $t\in \left[ 0,\psi \left( 0\right) %
\right] ,$ \\
$0,$ & for & $t\geq \psi \left( 0\right) .$%
\end{tabular}%
\ \right.
\end{equation*}%
When $\psi (0)=\infty ,$ the generator $\psi $ and $C$ are said to be
\textit{strict} and therefore $\psi ^{\lbrack -1]}=\psi ^{-1}.$ All notions
of positive dependence that appeared in the literature, including the
weakest one of PQD as defined by \cite{LEH66}, require the generator to be
strict. \medskip

\noindent Archimedean copulas are widely used in applications due to their
simple form, a variety of dependence structures and other \textquotedblleft
nice\textquotedblright\ properties. For example, in the Actuarial field: the
idea arose indirectly in \cite{Cla78} and was developed in \cite{Oak82},
\cite{CooJo81}. A survey of Actuarial applications is in \cite{FrVa98}%
.\medskip

\noindent For an Archimedean copula, the Kendall's tau can be evaluated
directly from the generator of the copula, as shown in \cite{Gen86}%
\begin{equation}
\tau =4{\displaystyle\int_{0}^{1}}\frac{\psi \left( u\right) }{\psi ^{\prime
}\left( u\right) }du+1.  \label{tau}
\end{equation}%
where $\psi ^{\prime }\left( u\right) $ exists a.e., since the generator is
convex. This is another \textquotedblleft nice\textquotedblright\ feature of
Archimedean copulas. As for tail dependency, as shown in \cite[page 105]%
{Joe97} the coefficient of upper tail dependency is
\begin{equation*}
\lambda _{U}=2-2\lim_{s\rightarrow 0^{+}}\frac{\psi ^{-1\prime }\left(
2s\right) }{\psi ^{-1\prime }\left( s\right) }
\end{equation*}%
and the coefficient of lower tail dependency is%
\begin{equation*}
\lambda _{L}=2\lim_{s\rightarrow +\infty }\frac{\psi ^{-1\prime }\left(
2s\right) }{\psi ^{-1\prime }\left( s\right) }.
\end{equation*}%
A collection of twenty-two one-parameter families of Archimedean copulas can
be found in Table 4.1 of \cite{Nels06}.\medskip

\noindent Notice that in the case of Archimedean copula the copula
conditional tail expectation has not an explicit formula, so we give by the
following Proposition the expression of CCTE in terms of generator.

\begin{proposition}
\label{JTUARCH}Let $C$ be an Archimedean copula absolutely continuous with
generator $\psi ,$ the CCTE of the target risk in terms of generator with
respect to the bivariate thresholds $(s,t),$ $0<s,t<1,$ is given by{}%
\begin{equation}
\mathbb{CCTE}_{X_{1}}\left( s;t\right) =\frac{1}{\overline{C}\left(
1-s,1-t\right) }\left( \left( 1-s\right) \mathbb{CTE}_{X_{1}}\left( s\right)
-\int_{s}^{1}\frac{\psi ^{\prime }(u)F_{X_{1}}^{-1}\left( u\right) }{\psi
^{\prime }\left( C\left( u,t\right) \right) }du\right) .  \label{archi}
\end{equation}
\end{proposition}

\noindent Note that in practice we can easily fit copula-based models with
the maximum likelihood method or with estimate the dependence parameter by
the relationship between Kendall's tau of the data and the generator of the
Archimedean copula given in (\ref{tau}) under the specified copula
model.\medskip

\noindent In the following section we give same examples to explain how to
calculate and compare the CCTE with other risk measure such VaR and CTE.

\subsubsection{\textbf{CCTE via Gumbel Copula}}

\noindent The Gumbel family has been introduced by \cite{Gu60}. Since it has
been discussed in \cite{Hou86}, it is also known as the Gumbel-Hougaard
family. The Gumbel copula is an asymmetric Archimedean copula. This copula
is given by%
\begin{equation*}
C_{\theta }^{G}\left( u,v\right) =\exp \left\{ -\left[ \left( -\ln u\right)
^{\theta }+\left( -\ln v\right) ^{\theta }\right] ^{1/\theta }\right\} ,
\end{equation*}%
its generator is%
\begin{equation*}
\psi _{\theta }\left( t\right) =\left( -\ln t\right) ^{\theta }.
\end{equation*}%
The dependence parameter is restricted to the interval $[1,\infty ).$ It
follows that the Gumbel family can represent independence and
\textquotedblleft positive\textquotedblright\ dependence only, since the
lower and upper bound for its parameter correspond to the product copula and
the upper Fr\'{e}chet bound. The Gumbel copula families is often used for
modeling heavy dependencies in right tail. It exhibits strong upper tail
dependence $\lambda _{U}=2-2^{1/\theta }$ and relatively weak lower tail
dependence $\lambda _{L}=0.$ If outcomes are known to be strongly correlated
at high values but less correlated at low values, then the Gumbel copula
will be an appropriate choice.\medskip

\noindent We give the CCTE of rv's $X_{i},$ $i=1,2,3$ in terms of Gumbel
copula by%
\begin{align}
\mathbb{CCTE}_{X_{i}}\left( s;t\right) & =\frac{1}{\overline{C}_{\theta
_{i}}^{G}\left( 1-s,1-t\right) }\left( \frac{\alpha \left( 1-s\right)
^{1-1/\alpha }}{\alpha -1}\right.  \notag \\
& \ \ \ \ \ \ \ \ \ \ \left. -\int_{s}^{1}u^{-1}\left( 1-u\right)
^{-1/\alpha }\left( -\ln u\right) ^{\theta _{i}-1}C_{\theta _{i}}^{G}\left(
u,t\right) \left( -\ln \left( C_{\theta _{i}}^{G}\left( u,t\right) \right)
\right) ^{1-\theta _{i}}du\right) ,  \label{Gumccte}
\end{align}%
where $\overline{C}_{\theta _{i}}^{G}\left( s,t\right) =s+t-1+C_{\theta
_{i}}^{G}(1-s,1-t).$\medskip

\noindent The CTE's and VaR's of $X_{i}$ is the same and it's given
respectively by (\ref{cteXi}) and (\ref{VaR}), for $i=1,2,3.$ \medskip

\noindent By the relationship between Kendall's tau $\tau $ and the Gumbel
copula parameter $\theta _{i}$ given by:%
\begin{equation*}
\tau =\left( \theta _{i}-1\right) /\theta _{i},
\end{equation*}%
we select the values of $\theta _{i}$ corresponding respectively to a weak,
a moderate and a strong positive association witch summarized in Table \ref{TAB2 Gumbel}.

\begin{table}[h] \centering%
\begin{tabular}{ccc}
$\lambda _{U}$ & $\theta _{i}$ & $\tau $ \\ \hline\hline
$0.013$ & $1.01$ & $0.009$ \\
$0.585$ & $2$ & $0.500$ \\
$0.928$ & $10$ & $0.900$ \\ \hline\hline
\multicolumn{1}{l}{} & \multicolumn{1}{l}{} & \multicolumn{1}{l}{}%
\end{tabular}%
\caption{Upper tail, Kendall's tau and Gumbel copula parameters used in calculate of risk measures.}\label{TAB2 Gumbel}%
\end{table}%

\begin{figure}[h]
\begin{center}
\includegraphics[height=3.5243in,width=4.754in]{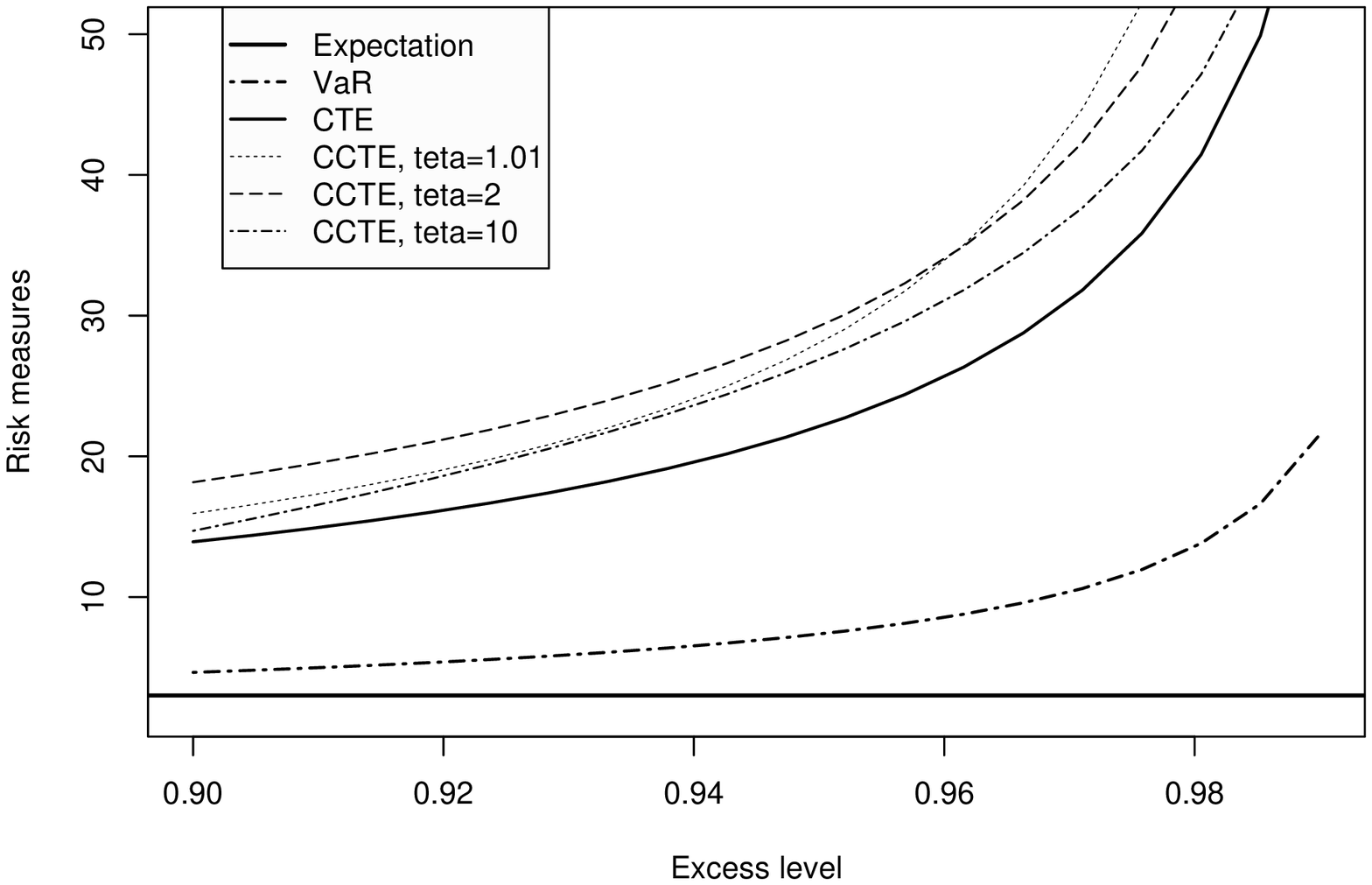}
\end{center}
\caption{$\mathbb{CCTE}$, $\mathbb{CTE}$ and $VaR$ risks measures of PQD
pareto (1.5) rv's with Gumbel copula and $0.9\leq s=t\leq 0.99.$}
\label{FigGumbel}
\end{figure}

\noindent Table \ref{TAB3}\ and Figure \ref{FigGumbel} shows that the loss $%
X_{1}$ is considerably riskier than $X_{2}$ and $X_{3},$ in the 6th column
of Table \ref{TAB3}, the relative difference between $112.1868$ and $69.6017$
is about $61.184\%.$\medskip

\begin{table}[h] \centering%
\begin{tabular}{rccccc}
${\small s}$ & \multicolumn{1}{|c}{${\small 0.9000}$} & ${\small 0.9225}$ & $%
{\small 0.9450}$ & ${\small 0.9675}$ & ${\small 0.9900}$ \\ \hline
${\small VaR}_{X_{i}}\left( s\right) $ & \multicolumn{1}{|c}{${\small 4.641}$%
} & ${\small 5.501}$ & ${\small 6.914}$ & ${\small 9.819}$ & ${\small 21.544}
$ \\
$\mathbb{CTE}_{X_{i}}\left( s\right) $ & \multicolumn{1}{|c}{${\small 13.924}
$} & ${\small 16.503}$ & ${\small 20.743}$ & ${\small 29.457}$ & ${\small %
64.633}$ \\ \hline
\multicolumn{1}{l}{${\small t}$} & \multicolumn{5}{|c}{$\mathbb{CCTE}%
_{X_{1}}\left( s,t\right) ,\ \ {\small \theta =1.01}$} \\ \hline
\multicolumn{1}{l}{${\small 0.9000}$} & \multicolumn{1}{|c}{${\small 15.937}$%
} & ${\small 16.485}$ & ${\small 17.410}$ & ${\small 19.365}$ & ${\small %
25.007}$ \\
\multicolumn{1}{l}{${\small 0.9225}$} & \multicolumn{1}{|c}{${\small 18.879}$%
} & ${\small 19.528}$ & ${\small 20.625}$ & ${\small 22.948}$ & ${\small %
33.690}$ \\
\multicolumn{1}{l}{${\small 0.9450}$} & \multicolumn{1}{|c}{${\small 23.699}$%
} & ${\small 24.507}$ & ${\small 25.873}$ & ${\small 28.760}$ & ${\small %
40.588}$ \\
\multicolumn{1}{l}{${\small 0.9675}$} & \multicolumn{1}{|c}{${\small 33.556}$%
} & ${\small 34.667}$ & ${\small 36.534}$ & ${\small 40.454}$ & ${\small %
56.275}$ \\
\multicolumn{1}{l}{${\small 0.9900}$} & \multicolumn{1}{|c}{${\small 72.992}$%
} & ${\small 75.133}$ & ${\small 78.645}$ & ${\small 85.726}$ & ${\small %
112.18}$ \\ \hline\hline
\multicolumn{1}{l}{${\small t}$} & \multicolumn{5}{|c}{$\mathbb{CCTE}%
_{X_{2}}\left( s,t\right) ,\ {\small \ \theta =2}$} \\ \hline
\multicolumn{1}{l}{${\small 0.9000}$} & \multicolumn{1}{|c}{${\small 18.158}$%
} & ${\small 19.769}$ & ${\small 22.691}$ & ${\small 28.950}$ & ${\small %
52.929}$ \\
\multicolumn{1}{l}{${\small 0.9225}$} & \multicolumn{1}{|c}{${\small 20.209}$%
} & ${\small 21.653}$ & ${\small 24.338}$ & ${\small 30.607}$ & ${\small %
53.742}$ \\
\multicolumn{1}{l}{${\small 0.9450}$} & \multicolumn{1}{|c}{${\small 23.842}$%
} & ${\small 25.059}$ & ${\small 27.383}$ & ${\small 33.070}$ & ${\small %
55.276}$ \\
\multicolumn{1}{l}{${\small 0.9675}$} & \multicolumn{1}{|c}{${\small 31.849}$%
} & ${\small 32.766}$ & ${\small 34.543}$ & ${\small 39.128}$ & ${\small %
59.207}$ \\
\multicolumn{1}{l}{${\small 0.9900}$} & \multicolumn{1}{|c}{${\small 66.087}$%
} & ${\small 66.606}$ & ${\small 67.583}$ & ${\small 70.074}$ & ${\small %
86.385}$ \\ \hline\hline
\multicolumn{1}{l}{${\small t}$} & \multicolumn{5}{|c}{$\mathbb{CCTE}%
_{X_{3}}\left( s,t\right) ,\ \ {\small \theta =10}$} \\ \hline
\multicolumn{1}{l}{${\small 0.9000}$} & \multicolumn{1}{|c}{${\small 13.765}$%
} & ${\small 16.694}$ & ${\small 23.338}$ & ${\small 39.483}$ & ${\small %
128.31}$ \\
\multicolumn{1}{l}{${\small 0.9225}$} & \multicolumn{1}{|c}{${\small 15.612}$%
} & ${\small 16.626}$ & ${\small 21.902}$ & ${\small 36.924}$ & ${\small %
120.00}$ \\
\multicolumn{1}{l}{${\small 0.9450}$} & \multicolumn{1}{|c}{${\small 19.378}$%
} & ${\small 19.446}$ & ${\small 20.821}$ & ${\small 32.807}$ & ${\small %
106.54}$ \\
\multicolumn{1}{l}{${\small 0.9675}$} & \multicolumn{1}{|c}{${\small 29.457}$%
} & ${\small 29.458}$ & ${\small 29.480}$ & ${\small 31.692}$ & ${\small %
95.737}$ \\
\multicolumn{1}{l}{${\small 0.9900}$} & \multicolumn{1}{|c}{${\small 64.633}$%
} & ${\small 65.034}$ & ${\small 66.412}$ & ${\small 67.753}$ & ${\small %
69.601}$ \\ \hline\hline
\multicolumn{1}{l}{} &  &  &  &  &
\end{tabular}%
\caption{Risk measures of PQD pareto (1.5) rv's with Gumbel copula.}\label{TAB3}%
\end{table}%

\noindent By definition of our risk measurement, the risks concern the study
is necessarily comonotonic, then to have a good decision we must select a
copula model with upper tail dependence, we show in next example that the
dependence models with no upper tails dependence do not helps us to take a
decision.

\subsubsection{\textbf{CCTE via Clayton Copula}}

\noindent In the following example, we consider the bivariate Clayton copula
which is a member of the class of Archimedean copula, with the dependence
parameter $\theta$ in $\left. \left[ -1,\infty\right) \right\backslash
\left\{ 0\right\} $.\medskip

\noindent The Clayton family was first proposed by \cite{Cla78} and studied
by \citeauthor{Oak82}, (\citeyear{Oak82}; \citeyear{Oak86}), %
\citeauthor{COa84}, (\citeyear{CooJo81}; \citeyear{CooJo86}). The Clayton
copula has been used to study correlated risks, it has the form%
\begin{equation}
C_{\theta }^{C}(u,v):=\left[ \max \left( u^{-\theta }+v^{-\theta
}-1,0\right) \right] ^{-1/\theta }.  \label{clayton}
\end{equation}%
For $\theta >0$ the copulas are strict and the copula expression simplifies
to
\begin{equation}
C_{\theta }^{C}(u,v)=\left( u^{-\theta }+v^{-\theta }-1\right) ^{-1/\theta }.
\label{claytoncop}
\end{equation}%
Asymmetric tail dependence is prevalent if the probability of joint extreme
negative realizations differs from that of joint extreme positive
realizations. it can be seen that the Clayton copula assigns a higher
probability to joint extreme negative events than to joint extreme positive
events. The Clayton copula is said to display lower tail dependence $\lambda
_{L}=2^{-1/\theta },$ while it displays zero upper tail dependence $\lambda
_{U}=0,$ for $\theta \geq 0.$ The converse can be said about the Gumbel
copula (displaying upper but zero lower tail dependence). The margins become
independent as $\theta $ approaches to zero, while for $\theta \rightarrow
\infty ,$ the Clayton copula arrives at the comonotonicity copula. For $%
\theta =-1$ we obtain the Fr\'{e}chet-Hoeffding lower bound and the copula
attains the Fr\'{e}chet upper bound as $\theta $ approaches to
infinity.\medskip

\noindent Clayton copula is the best suited for applications in which two
outcomes are likely to experience low values together, since the dependence
is strong in the lower tail and weak in the upper tail.\medskip

\noindent We take the same example as in the Subsection \ref{FGMsubsection},
we may now represents the joint df's $H_{i},$ $i=1,2,3,$\ respectively by
the Clayton copulas $C_{\theta _{i}}^{C}$ given in (\ref{claytoncop}) to
have an idea about the effects of lower tail dependence on our risk
measurement.\medskip

\noindent The relationship between Kendall's tau $\tau $ and the Clayton
copula is given by%
\begin{equation}
\tau =\theta _{i}/\left( \theta _{i}+2\right) ,  \label{tau Clay}
\end{equation}%
we select a different dependents parameters corresponding to several levels
of positive dependency summarized in Table \ref{TAB1 clayton} for a weak, a
moderate and a strong positive association, to calculate and compare the
CCTE's of $X_{i},i=1,2,3.$

\begin{table}[h] \centering%
\begin{tabular}{ccc}
$\lambda _{L}$ & $\theta _{i}$ & $\tau $ \\ \hline\hline
$0.250$ & $0.5$ & $0.200$ \\
$0.707$ & $2$ & $0.500$ \\
$0.943$ & $12$ & $0.857$ \\ \hline\hline
\multicolumn{1}{l}{} & \multicolumn{1}{l}{} & \multicolumn{1}{l}{}%
\end{tabular}%
\caption{Lower tail, Kendall's tau and Clayton copula parameters used in calculate of risk measures.}\label{TAB1 clayton}
\end{table}%
%EndExpansion

\begin{figure}[h]
\begin{center}
\includegraphics[height=3.5243in,width=4.754in]{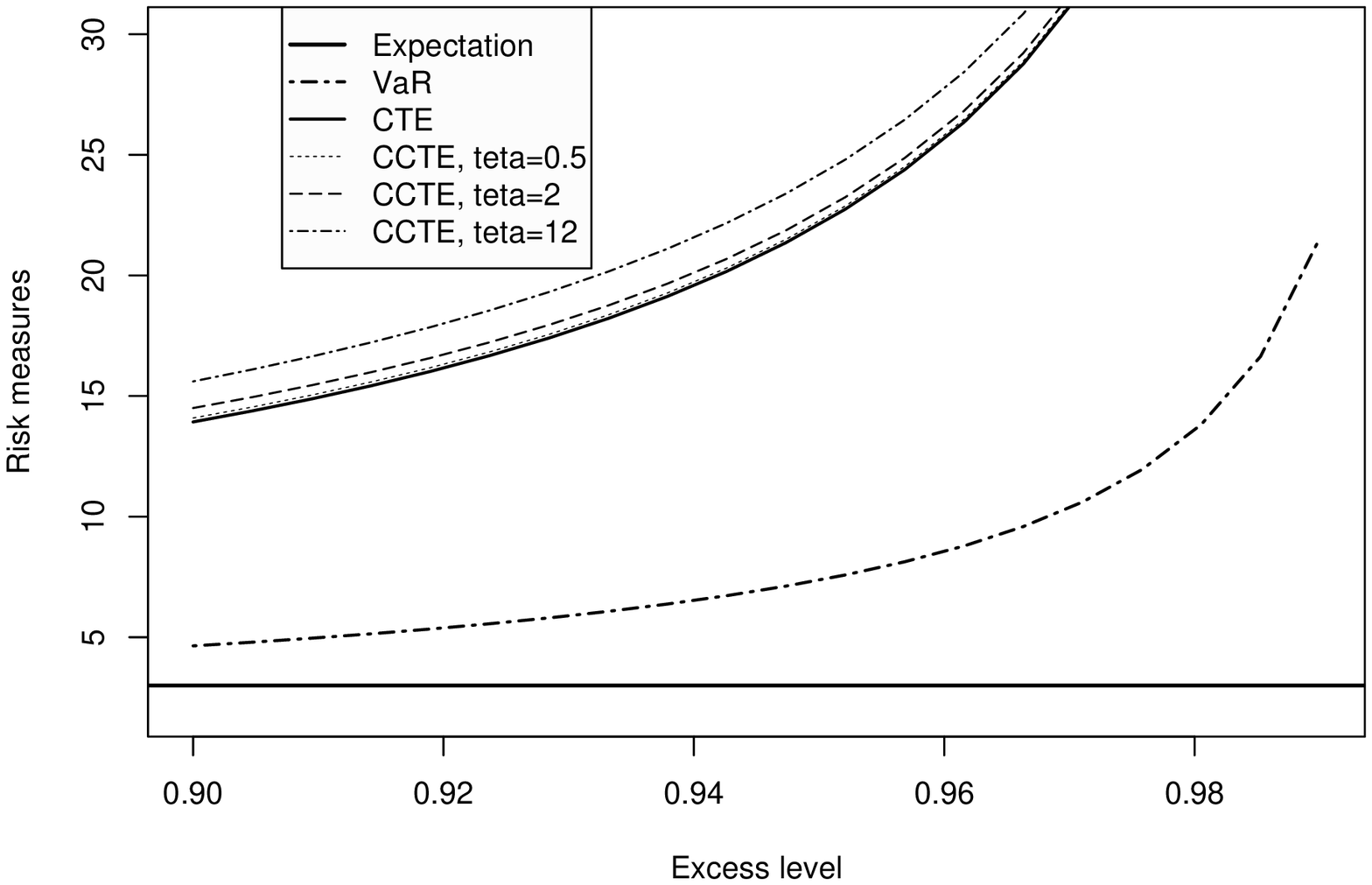}
\end{center}
\caption{$\mathbb{CCTE}$, $\mathbb{CTE}$ and $VaR$ risks measures of PQD
pareto (1.5) rv's with Clayton copula and $0.9\leq s=t\leq0.99.$}
\label{FigClayton}
\end{figure}

\noindent The CCTE of the rv's $X_{i}$ with respect to the bivariate
thresholds $(s,t)$ is given by%
\begin{equation}
\mathbb{CCTE}_{X_{i}}\left( s;t\right) =\frac{1}{\overline{C}_{\theta
_{i}}^{C}\left( 1-s,1-t\right) }\left( \frac{\alpha \left( 1-s\right)
^{-1/\alpha +1}}{\left( \alpha -1\right) }-{\displaystyle\int_{s}^{1}}\frac{%
\left( t^{-\theta _{i}}+u^{-\theta _{i}}-1\right) ^{-1-1/\theta _{i}}}{%
\left( 1-u\right) ^{1/\alpha }u^{\theta _{i}+1}}du\right) .  \label{Claccte}
\end{equation}%
The differences as reported in Table \ref{TAB2} and Figure \ref{FigClayton}
do not look very significant, in the 6th column of Table \ref{TAB2}, the
relative difference between $66.3802$ and $64.6330$ is only about $1.027\%.$
The differences is not found also when $t$ is small compared to $s,$ $%
\mathbb{CCTE}_{X_{1}}\left( 0.99;0.01\right) =64.6332$ and $\mathbb{CCTE}%
_{X_{3}}\left( 0.99;0.01\right) =64.6329$ the difference is about $1\%.$%
\medskip

\begin{table}[h] \centering
\begin{tabular}{rccccc}
$s$ & \multicolumn{1}{|c}{${\small 0.9000}$} & ${\small 0.9225}$ & ${\small %
0.9450}$ & ${\small 0.9675}$ & ${\small 0.9900}$ \\ \hline
${\small VaR}_{X_{i}}\left( s\right) $ & \multicolumn{1}{|c}{${\small 4.641}$%
} & ${\small 5.501}$ & ${\small 6.914}$ & ${\small 9.819}$ & ${\small 21.544}
$ \\
$\mathbb{CTE}_{X_{i}}\left( s\right) $ & \multicolumn{1}{|c}{${\small 13.924}
$} & ${\small 16.503}$ & ${\small 20.743}$ & ${\small 29.457}$ & ${\small %
64.633}$ \\ \hline
\multicolumn{1}{l}{$t$} & \multicolumn{5}{|c}{$\mathbb{CCTE}_{X_{1}}\left(
s,t\right) ,\ \ {\small \theta =0.5}$} \\ \hline
\multicolumn{1}{l}{${\small 0.9000}$} & \multicolumn{1}{|c}{${\small 14.088}$%
} & ${\small 14.092}$ & ${\small 14.096}$ & ${\small 14.101}$ & ${\small %
14.105}$ \\
\multicolumn{1}{l}{${\small 0.9225}$} & \multicolumn{1}{|c}{${\small 16.652}$%
} & ${\small 16.656}$ & ${\small 16.660}$ & ${\small 16.664}$ & ${\small %
16.667}$ \\
\multicolumn{1}{l}{${\small 0.9450}$} & \multicolumn{1}{|c}{${\small 20.874}$%
} & ${\small 20.878}$ & ${\small 20.881}$ & ${\small 20.884}$ & ${\small %
20.888}$ \\
\multicolumn{1}{l}{${\small 0.9675}$} & \multicolumn{1}{|c}{${\small 29.566}$%
} & ${\small 29.569}$ & ${\small 29.572}$ & ${\small 29.575}$ & ${\small %
29.577}$ \\
\multicolumn{1}{l}{${\small 0.9900}$} & \multicolumn{1}{|c}{${\small 64.706}$%
} & ${\small 64.707}$ & ${\small 64.709}$ & ${\small 64.711}$ & ${\small %
64.713}$ \\ \hline\hline
\multicolumn{1}{l}{$t$} & \multicolumn{5}{|c}{$\mathbb{CCTE}_{X_{2}}\left(
s,t\right) ,\ \ {\small \theta =2}$} \\ \hline
\multicolumn{1}{l}{${\small 0.9000}$} & \multicolumn{1}{|c}{${\small 14.500}$%
} & ${\small 14.536}$ & ${\small 14.572}$ & ${\small 14.610}$ & ${\small %
14.648}$ \\
\multicolumn{1}{l}{${\small 0.9225}$} & \multicolumn{1}{|c}{${\small 17.023}$%
} & ${\small 17.056}$ & ${\small 17.089}$ & ${\small 17.123}$ & ${\small %
17.159}$ \\
\multicolumn{1}{l}{${\small 0.9450}$} & \multicolumn{1}{|c}{${\small 21.199}$%
} & ${\small 21.227}$ & ${\small 21.257}$ & ${\small 21.288}$ & ${\small %
21.319}$ \\
\multicolumn{1}{l}{${\small 0.9675}$} & \multicolumn{1}{|c}{${\small 29.833}$%
} & ${\small 29.857}$ & ${\small 29.882}$ & ${\small 29.907}$ & ${\small %
29.934}$ \\
\multicolumn{1}{l}{${\small 0.9900}$} & \multicolumn{1}{|c}{${\small 64.882}$%
} & ${\small 64.898}$ & ${\small 64.915}$ & ${\small 64.932}$ & ${\small %
64.950}$ \\ \hline\hline
\multicolumn{1}{l}{$t$} & \multicolumn{5}{|c}{$\mathbb{CCTE}_{X_{3}}\left(
s,t\right) ,\ \ {\small \theta =12}$} \\ \hline
\multicolumn{1}{l}{${\small 0.9000}$} & \multicolumn{1}{|c}{${\small 15.605}$%
} & ${\small 16.118}$ & ${\small 16.743}$ & ${\small 17.494}$ & ${\small %
18.383}$ \\
\multicolumn{1}{l}{${\small 0.9225}$} & \multicolumn{1}{|c}{${\small 17.913}$%
} & ${\small 18.366}$ & ${\small 18.930}$ & ${\small 19.618}$ & ${\small %
20.447}$ \\
\multicolumn{1}{l}{${\small 0.9450}$} & \multicolumn{1}{|c}{${\small 21.888}$%
} & ${\small 22.274}$ & ${\small 22.762}$ & ${\small 23.371}$ & ${\small %
24.119}$ \\
\multicolumn{1}{l}{${\small 0.9675}$} & \multicolumn{1}{|c}{${\small 30.331}$%
} & ${\small 30.637}$ & ${\small 31.033}$ & ${\small 31.536}$ & ${\small %
32.169}$ \\
\multicolumn{1}{l}{${\small 0.9900}$} & \multicolumn{1}{|c}{${\small 65.169}$%
} & ${\small 65.363}$ & ${\small 65.619}$ & ${\small 65.951}$ & ${\small %
66.380}$ \\ \hline\hline
\multicolumn{1}{l}{} &  &  &  &  &
\end{tabular}%
\caption{Risk measures of dependent pareto (1.5) rv's with Clayton copula.}\label{TAB2}%
\end{table}%

\section{Application\label{Applicaltion}}

\noindent The relationship between the parameter of an Archimedean copula
and Kendall's tau has allowed us to calculate the value of this parameter
assuming a well precise Archimedean copula e.g., Gumbel copula. Once endowed
with the parameter value, we are able to compute any joint probability
between the stock indices.\medskip

\begin{figure}[h]
\begin{center}
\includegraphics[height=3.5218in,width=3.5309in]{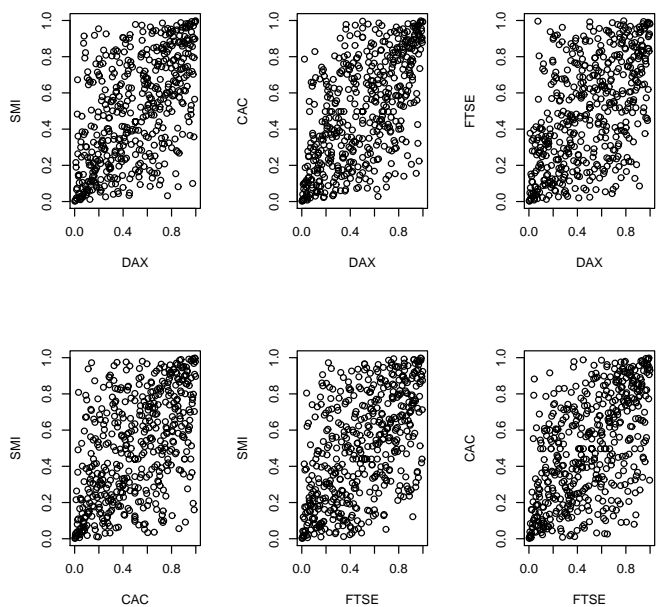}
\end{center}
\caption{Scatterplots of $500$ pseudo-observations drawn from a four
European stock indices returns.}
\label{Fig DSCF}
\end{figure}

\noindent For instance we analyzed $500$ observations from four European
stock indices return series calculated by $\log \left( X_{t+1}/X_{t}\right) $
for the period July 1991 to June 1993 (see, Figure \ref{Fig DSCF} ),
available in "QRM and datasets packages" of R software, it contains the
daily closing prices of major European stock indices: Germany DAX (Ibis),
Switzerland SMI, France CAC and UK FTSE. The data are sampled in business
time, i.e., weekends and holidays are omitted. Table \ref{TABS1} summaries
the Kendall's tau between the four Market Index returns.

\begin{table}[h] \centering
\begin{tabular}{ccccc}
\hline
\multicolumn{1}{|c}{\small Variable} & \multicolumn{1}{|c}{\small DAX} &
\multicolumn{1}{|c}{\small SMI} & \multicolumn{1}{|c}{\small CAC} &
\multicolumn{1}{|c|}{\small FTSE} \\ \hline\hline
\multicolumn{1}{|c}{\small DAX} & \multicolumn{1}{|c}{${\small 1}$} &
\multicolumn{1}{|c}{${\small 0.4052}$} & \multicolumn{1}{|c}{${\small 0.4374}
$} & \multicolumn{1}{|c|}{${\small 0.3706}$} \\ \hline
\multicolumn{1}{|c}{\small SMI} & \multicolumn{1}{|c}{${\small 0.4052}$} &
\multicolumn{1}{|c}{${\small 1}$} & \multicolumn{1}{|c}{${\small 0.3791}$} &
\multicolumn{1}{|c|}{${\small 0.3924}$} \\ \hline
\multicolumn{1}{|c}{\small CAC} & \multicolumn{1}{|c}{${\small 0.4374}$} &
\multicolumn{1}{|c}{${\small 0.3791}$} & \multicolumn{1}{|c}{${\small 1}$} &
\multicolumn{1}{|c|}{${\small 0.4076}$} \\ \hline
\multicolumn{1}{|c}{\small FTSE} & \multicolumn{1}{|c}{${\small 0.3706}$} &
\multicolumn{1}{|c}{${\small 0.3924}$} & \multicolumn{1}{|c}{${\small 0.4076}
$} & \multicolumn{1}{|c|}{${\small 1}$} \\ \hline\hline
&  &  &  &
\end{tabular}%
\caption{Kendall's tau matrix estimates from four European stock indices returns.}\label{TABS1}%
\end{table}%

\noindent By assuming that Gumbel copula represents our four dependence
structures, we obtain the fitted dependence parameters of the six bivariate
joint df's, presented in Table \ref{TABS2}.

\begin{table}[h] \centering%
\begin{tabular}{ccccc}
\hline
\multicolumn{1}{|c}{\small Variable} & \multicolumn{1}{|c}{\small DAX} &
\multicolumn{1}{|c}{\small SMI} & \multicolumn{1}{|c}{\small CAC} &
\multicolumn{1}{|c||}{\small FTSE} \\ \hline\hline
\multicolumn{1}{|c}{\small DAX} & \multicolumn{1}{|c}{${\small \infty }$} &
\multicolumn{1}{|c}{${\small 1.6815}$} & \multicolumn{1}{|c}{${\small 1.7777}
$} & \multicolumn{1}{|c||}{${\small 1.5888}$} \\ \hline
\multicolumn{1}{|c}{\small SMI} & \multicolumn{1}{|c}{${\small 1.6815}$} &
\multicolumn{1}{|c}{${\small \infty }$} & \multicolumn{1}{|c}{${\small 1.6106%
}$} & \multicolumn{1}{|c||}{${\small 1.6459}$} \\ \hline
\multicolumn{1}{|c}{\small CAC} & \multicolumn{1}{|c}{${\small 1.7777}$} &
\multicolumn{1}{|c}{${\small 1.6106}$} & \multicolumn{1}{|c}{${\small \infty
}$} & \multicolumn{1}{|c||}{${\small 1.6880}$} \\ \hline
\multicolumn{1}{|c}{\small FTSE} & \multicolumn{1}{|c}{${\small 1.5888}$} &
\multicolumn{1}{|c}{${\small 1.6459}$} & \multicolumn{1}{|c}{${\small 1.6880}
$} & \multicolumn{1}{|c||}{${\small \infty }$} \\ \hline\hline
&  &  &  &
\end{tabular}%
\caption{Fitted copula parameter correspoding to Kendall's tau, Gumbel copula.}\label{TABS2}%
\end{table}%

\noindent The $\alpha $-stable distribution offers a reasonable improvement
to the alternative distributions, each stable distribution $S_{\alpha
}(\sigma ;\beta ;\mu )$ has the stability index $\alpha $ that can be
treated as the main parameter, when we make an investment decision, skewness
parameter $\beta $, in the range $[-1,1]$, scale parameter $\sigma $ and
shift parameter $\mu .$ In models that use financial data, it is generally
assumed that $\alpha \in (1,2].$ By using the "fBasics" package in R
software, based on the maximum likelihood estimators to fit the parameters
of a df's of the four Market Index returns, the results are summarized in
Table \ref{TABS3}.

%TCIMACRO{\TeXButton{B}{\begin{table}[h] \centering}}%
%BeginExpansion
\begin{table}[h] \centering%
%EndExpansion
\begin{tabular}{ccccc}
& {\small DAX} & {\small SMI} & {\small CAC} & {\small FTSE} \\ \hline\hline
${\small \alpha }$ & \multicolumn{1}{r}{${\small 1.6420}$} &
\multicolumn{1}{r}{${\small 1.8480}$} & \multicolumn{1}{r}{${\small 1.6930}$}
& \multicolumn{1}{r}{${\small 1.8740}$} \\
${\small \beta }$ & \multicolumn{1}{r}{${\small 0.1470}$} &
\multicolumn{1}{r}{${\small 0.1100}$} & \multicolumn{1}{r}{${\small 0.0380}$}
& \multicolumn{1}{r}{${\small 0.9500}$} \\
${\small \sigma }$ & \multicolumn{1}{r}{${\small 0.0046}$} &
\multicolumn{1}{r}{${\small 0.0045}$} & \multicolumn{1}{r}{${\small 0.0006}$}
& \multicolumn{1}{r}{${\small 0.0053}$} \\
${\small \mu }$ & \multicolumn{1}{r}{${\small -0.0001}$} &
\multicolumn{1}{r}{${\small 0.0006}$} & \multicolumn{1}{r}{${\small -0.0001}$%
} & \multicolumn{1}{r}{${\small -0.0005}$} \\ \hline\hline
&  &  &  &
\end{tabular}%
\caption{Maximum likelihood fit of four-parameters stable distribution to
four European stock indices retuns data.}\label{TABS3}%
%TCIMACRO{\TeXButton{E}{\end{table}}}%
%BeginExpansion
\end{table}%
%EndExpansion

\noindent The $\alpha $-stable distribution has Pareto-type tails, it's like
a power function, i.e., $F$ is regularly varying (at infinity) with index $%
\left( -\alpha \right) ,$ meaning that $\overline{F}\left( x\right)
=x^{-\alpha }L\left( x\right) $ as $x$ becomes large, where $L>0$ is a
slowly varying function, which can be interpreted as slower than any power
function (see, \citeauthor{Res87}; \citeyear{Res87} and \citeauthor{Sen76}; %
\citeyear{Sen76} for a technical treatment of regular variation). By using
the Equations (\ref{Gumccte}) for the Gumbel copula fitting, we calculate
for a fixed levels $s=t=0.95$ the CCTE's risk measures for the all cases,
the results are summarized in Table \ref{TABS4}.

%TCIMACRO{\TeXButton{B}{\begin{table}[h] \centering}}%
%BeginExpansion
\begin{table}[h] \centering%
%EndExpansion
\begin{tabular}{ccccc}
\hline
\multicolumn{1}{|c}{\small Variable} & \multicolumn{1}{|c}{\small DAX} &
\multicolumn{1}{|c}{\small SMI} & \multicolumn{1}{|c}{\small CAC} &
\multicolumn{1}{|c||}{\small FTSE} \\ \hline\hline
\multicolumn{1}{|c}{\small DAX} & \multicolumn{1}{|c}{${\small -}$} &
\multicolumn{1}{|c}{${\small 21.5009}$} & \multicolumn{1}{|c}{${\small %
21.0786}$} & \multicolumn{1}{|c||}{${\small 21.9731}$} \\ \hline
\multicolumn{1}{|c}{\small SMI} & \multicolumn{1}{|c}{${\small 14.2812}$} &
\multicolumn{1}{|c}{${\small -}$} & \multicolumn{1}{|c}{${\small 14.4703}$}
& \multicolumn{1}{|c||}{${\small 14.3737}$} \\ \hline
\multicolumn{1}{|c}{\small CAC} & \multicolumn{1}{|c}{${\small 18.8362}$} &
\multicolumn{1}{|c}{${\small 19.4915}$} & \multicolumn{1}{|c}{${\small -}$}
& \multicolumn{1}{|c||}{${\small 19.1671}$} \\ \hline
\multicolumn{1}{|c}{\small FTSE} & \multicolumn{1}{|c}{${\small 13.9075}$} &
\multicolumn{1}{|c}{${\small 13.7593}$} & \multicolumn{1}{|c}{${\small %
13.6576}$} & \multicolumn{1}{|c||}{${\small -}$} \\ \hline\hline
&  &  &  &
\end{tabular}%
\caption{CCTE's Risk measures for $s=0.99$ and $t=0.99$ with Gumbel copula
(left panel) and Clayton copula (right panel).}\label{TABS4}%
%TCIMACRO{\TeXButton{E}{\end{table}}}%
%BeginExpansion
\end{table}%
%EndExpansion

\noindent The smallest values in Table \ref{TABS4} gives the lowest risk.
So, the less risky couples $(X,Y)$ are: (DAX, CAC), (SMI, DAX), (CAC, DAX)
and (FTSE, CAC), where $X$ is the target risk and $Y$ is the associated risk.

\section{\label{Conclu}Conclusion notes}

\noindent This paper discussed a new risk measure called copula conditional
tail expectation. This measure aid to understanding the relationships among
multivariate assets and to help us significantly about how best to position
our investments and improve our financial risk protection.\medskip

\noindent Tables \ref{TAB3} show that the copula conditional tail
expectation measure become smaller as the dependency increase. However, CTE
and VaR are neither increasing nor decreasing as the correlation increase.
Therefore, the dependency information helps us to minimize the risk.\medskip

\noindent \textbf{Acknowledgements. }The author is indebted to an anonymous
referee for their careful reading and suggestions for improvements.

\section{Appendix}

\begin{proof}[\textbf{Proof of Proposition \protect\ref{Prop1}}]
By conditional probability is easily to obtain
\begin{equation*}
\mathbb{P}\left( \left. X_{1}\leq x\right\vert X_{1}>VaR_{X_{1}}\left(
s\right) ,X_{2}>VaR_{X_{2}}\left( t\right) \right) =\frac{\mathbb{P}\left(
x\geq X_{1}>VaR_{X_{1}}\left( s\right) ,X_{2}>VaR_{X_{2}}\left( t\right)
\right) }{\mathbb{P}\left( X_{1}>VaR_{X_{1}}\left( s\right)
,X_{2}>VaR_{X_{2}}\left( t\right) \right) }\medskip
\end{equation*}%
On the other hand, we have%
\begin{eqnarray*}
&&\mathbb{P}\left( x\geq X_{1}>VaR_{X_{1}}\left( s\right)
,X_{2}>VaR_{X_{2}}\left( t\right) \right) \\
&=&1-\mathbb{P}\left( F_{1}^{-1}\left( x\right) <F_{1}\left( X_{1}\right)
\leq s\right) -\mathbb{P}\left( F_{2}\left( X_{2}\right) \leq t\right) +%
\mathbb{P}\left( F_{1}^{-1}\left( x\right) <F_{1}\left( X_{1}\right) \leq
s,F_{2}\left( X_{2}\right) \leq t\right) ,
\end{eqnarray*}%
\begin{eqnarray*}
&&\mathbb{P}\left( X_{1}>VaR_{X_{1}}\left( s\right) ,X_{2}>VaR_{X_{2}}\left(
t\right) \right) \  \\
&=&1-\mathbb{P}\left( F_{1}\left( X_{1}\right) \leq s\right) -\mathbb{P}%
\left( F_{2}\left( X_{2}\right) \leq t\right) +\mathbb{P}\left( F_{1}\left(
X_{1}\right) \leq s,F_{2}\left( X_{2}\right) \leq t\right) \\
&=&1-s-t+C\left( s,t\right) \\
&=&\overline{C}\left( 1-s,1-t\right) ,
\end{eqnarray*}%
\begin{equation*}
\mathbb{P}\left( F_{1}^{-1}\left( x\right) <F_{1}\left( X_{1}\right) \leq
s\right) =s-VaR_{X_{1}}\left( x\right)
\end{equation*}%
and%
\begin{equation*}
\mathbb{P}\left( F_{1}^{-1}\left( x\right) <F_{1}\left( X_{1}\right) \leq
s,F_{2}\left( X_{2}\right) \leq t\right) =C\left( s,t\right)
-C(VaR_{X_{1}}(x),t).
\end{equation*}%
Then%
\begin{equation*}
\mathbb{P}\left( \left. X_{1}\leq x\right\vert X_{1}>VaR_{X_{1}}\left(
s\right) ,X_{2}>VaR_{X_{2}}\left( t\right) \right) =1+\frac{%
VaR_{X_{1}}\left( x\right) -C(VaR_{X_{1}}(x),t)}{\overline{C}\left(
1-s,1-t\right) }
\end{equation*}%
Then the CCTE is given by%
\begin{eqnarray*}
\mathbb{CCTE}_{X_{1}}\left( s,t\right) &=&\int_{X_{1}>VaR_{X_{1}}\left(
s\right) }xd\mathbb{P}\left( \left. X_{1}\leq x\right\vert
X_{1}>VaR_{X_{1}}\left( s\right) ,X_{2}>VaR_{X_{2}}\left( t\right) \right) \\
&=&\frac{1}{\overline{C}\left( 1-s,1-t\right) }\int_{VaR_{X_{1}}\left(
s\right) }^{\infty }xd\left( VaR_{X_{1}}\left( x\right)
-C(VaR_{X_{1}}(x),t)\right) . \\
&=&\frac{1}{\overline{C}\left( 1-s,1-t\right) }\int_{s}^{1}F_{X_{1}}^{-1}%
\left( u\right) d\left( u-C(u,t)\right) \\
&=&\frac{1}{\overline{C}\left( 1-s,1-t\right) }\left(
\int_{s}^{1}F_{X_{1}}^{-1}\left( u\right)
du-\int_{s}^{1}F_{X_{1}}^{-1}\left( u\right) dC(u,t)\right)
\end{eqnarray*}%
This close the proof of Proposition \ref{Prop1}.\bigskip
\end{proof}

\begin{proof}[\textbf{Proof of Proposition \protect\ref{JTUARCH}}]
Let's denote by
\begin{equation*}
C_{u}\left( u,v\right) :=\frac{\partial C\left( u,v\right) }{\partial u}
\end{equation*}%
then by (\ref{CCTE}), we have%
\begin{equation*}
\mathbb{CCTE}_{X_{1}}\left( s;t\right) =\frac{1}{\overline{C}\left(
1-s,1-t\right) }\left( \int_{s}^{1}F_{X_{1}}^{-1}\left( u\right)
du-\int_{s}^{1}F_{X_{1}}^{-1}\left( u\right) C_{u}(u,t)du\right) .
\end{equation*}%
So, $C$ is Archimedean copula, then%
\begin{equation}
C_{u}\left( u,v\right) =\frac{\psi ^{\prime }(u)}{\psi ^{\prime }\left(
C\left( u,v\right) \right) },  \label{archi2}
\end{equation}%
Finely, we get (\ref{archi}) by substitution of%
\begin{equation*}
\int_{s}^{1}F_{X_{1}}^{-1}\left( u\right) du=\left( 1-s\right) \mathbb{CTE}%
_{X_{1}}\left( s\right)
\end{equation*}%
and (\ref{archi2}) in (\ref{CCTE}).
\end{proof}

\end{document}